\documentclass[12pt]{amsart}
\usepackage{amsmath}
\usepackage{verbatim}
\usepackage{epsfig}
\usepackage{graphics}
\usepackage{latexsym}
\usepackage{amscd}

\newcommand{\bd}{\partial}
\newcommand{\hm}{ho\-me\-o\-mor\-phic}
\newcommand{\R}{\ensuremath{\mathbf{R}}}
\newcommand{\RR}{\ensuremath{\mathbf{R}^2}}
\newcommand{\RRR}{\ensuremath{\mathbf{R}^3}}
\newcommand{\sbs}{\subseteq}
\newcommand{\irr}{ir\-re\-du\-ci\-ble}
\newcommand{\inc}{in\-com\-press\-i\-ble}
\newcommand{\birr}{$\bd$-\irr}
\newcommand{\rirr}{\RR-\irr}
\newcommand{\pl}{parallel}
\newcommand{\bpl}{$\bd$-\pl}
\newcommand{\ei}{end \irr}

\newcommand{\tm}{3-manifold}

\newcommand{\p}{^{\prime}}
\newcommand{\pp}{^{\prime\prime}}
\newcommand{\er}{end reduction}
\newcommand{\ra}{\rightarrow}
\newcommand{\ns}{\emptyset}
\newcommand{\inte}{int \,}
\newcommand{\ga}{\ensuremath{\gamma}}
\newcommand{\be}{\ensuremath{\beta}}
\newcommand{\al}{\ensuremath{\alpha}}
\newcommand{\ze}{\ensuremath{\zeta}}

\newcommand{\ta}{\ensuremath{\tau}}

\newcommand{\8}{\ensuremath{\infty}}
\newcommand{\n}{^{-1}}
\newcommand{\Inte}{Int \, }
\newcommand{\halfspace}{\ensuremath{\RR \times [0,\8)}}

\newcommand{\ep}{\ensuremath{\varepsilon}}
\newcommand{\tht}{\ensuremath{\theta}}
\newcommand{\om}{\ensuremath{\omega}}
\newcommand{\ka}{\ensuremath{\kappa}}
\newcommand{\pe}{poly-excellent}
\newcommand{\ps}{poly-superb}
\newcommand{\de}{\ensuremath{\delta}}

\newcommand{\la}{\ensuremath{\lambda}}

\newcommand{\SW}{\ensuremath{\mathcal{S}(W)}}
\newcommand{\De}{\ensuremath{\Delta}}
\newcommand{\Dep}{\ensuremath{\De\p}}
\newcommand{\Depp}{\ensuremath{\De\pp}}
\newcommand{\Ga}{\ensuremath{\Gamma}}
\newcommand{\anan}{anannular}
\newcommand{\binc}{$\bd$-\inc}
\newcommand{\qe}{quasi-exhaustion}
\newcommand{\FF}{\ensuremath{\mathcal{F}}}
\newcommand{\EE}{\ensuremath{\mathcal{E}}}
\newcommand{\BB}{\ensuremath{\mathcal{B}}}
\newcommand{\PP}{\ensuremath{\mathcal{P}}}
\newcommand{\whEE}{\ensuremath{\widehat{\EE}}}
\newcommand{\whBB}{\ensuremath{\widehat{\BB}}}
\newcommand{\whPP}{\ensuremath{\widehat{\PP}}}
\newcommand{\whP}{\ensuremath{\widehat{P}}}

\newcommand{\whC}{\ensuremath{\widehat{C}}}
\newcommand{\whB}{\ensuremath{\widehat{B}}}
\newcommand{\whE}{\ensuremath{\widehat{E}}}
\newcommand{\DD}{\ensuremath{\mathcal{D}}}
\newcommand{\QQ}{\ensuremath{\mathcal{Q}}}
\newcommand{\SSS}{\ensuremath{\mathcal{S}}}
\newcommand{\AAA}{\ensuremath{\mathcal{A}}}
\newcommand{\YY}{\ensuremath{\mathcal{Y}}}
\newcommand{\NN}{\ensuremath{\mathcal{N}}}
\newcommand{\mgp}{minimal general position}
\newcommand{\sh}{^{\#}}

\title[The complex of end reductions]
{The complex of end reductions \\ of a contractible open 
3-manifold: \\ constructing 1-dimensional examples}
\author{Robert Myers} 
\address{Department of Mathematics, Oklahoma State University, Stillwater, 
OK 74078}
\email{myersr@math.okstate.edu}
\dedicatory{Dedicated to Fico Gonz\'{a}lez-Acu\~{n}a in honor of his 
60th birthday}
\thanks{This research was partially supported by NSF Grant DMS-0072429.}
\renewcommand{\labelenumi}{(\theenumi)}
\newtheorem{thm}{Theorem}[section]
\newtheorem{prop}[thm]{Proposition}
\newtheorem{cor}[thm]{Corollary}
\newtheorem{lem}[thm]{Lemma}

\begin{document}

\begin{abstract} 
Given an \irr, contractible, open \tm\ $W$ which is not \hm\ to \RRR, 
there is an associated simplicial complex $\mathcal{S}(W)$, the complex 
of \er s of $W$. Whenever $W$ covers a \tm\ $M$ one has that $\pi_1(M)$ 
is isomorphic to a subgroup of the group $Aut(\mathcal{S}(W))$ of 
simplicial automorphisms of $\mathcal{S}(W)$. 

In this paper we give a new method for constructing examples 
$W$ with $\mathcal{S}(W)$ isomorphic to a triangulation of \R. It follows 
that any \tm\ $M$ covered by $W$ must have $\pi_1(M)$ infinite 
cyclic. We also give a complete isotopy classification of the end 
reductions of these $W$. 
\end{abstract}

\maketitle

%Section 1
\section{Introduction}

A \textit{Whitehead manifold} $W$ is an \irr, contractible, open \tm\ 
which is not \hm\ to \RRR. Given a compact \tm\ $J$ in $W$ which is not 
contained in a 3-ball in $W$ Brin and Thickstun \cite{BT} defined a certain 
open submanifold $V$ of $W$ called an \er\ of $W$ at $J$. End reductions 
are rather nicely behaved but badly embedded manifolds which have certain 
interesting engulfing and homotopy theoretic properties and are unique up to 
isotopy with respect to these properties.

In \cite{My Endcov} the author showed how to associate to the set of 
isotopy classes of end reductions of $W$ a certain abstract simplicial 
complex $\mathcal{S}(W)$ with the following properties. Every 
self-homeomorphism of $W$ induces an automorphism of $\mathcal{S}(W)$. 
Whenever $W$ is a non-trivial covering space of a \tm\ $M$ each 
non-trivial element of the group $\pi_1(M)$ of covering 
translations acts without fixed points on \SW. Thus information 
about \SW\ gives information about what \tm s $W$ can cover.  

This complex seems particularly useful when $W$ is \rirr, i.e. 
when $W$ contains no ``non-trivial'' planes. In \cite{My Endcov} the author 
considered an uncountable collection of \rirr\ Whitehead manifolds 
which are modifications of an example due to Scott and Tucker \cite{ST}. 
He showed that each of these manifolds has \SW\ isomorphic to a 
triangulation of the real line. It follows that each \tm\ which 
is non-trivially covered by one of these \tm s must have infinite 
cyclic fundamental group, and in fact there are uncountably many which 
do cover such manifolds. 

These ``modified Scott-Tucker manifolds'' are easy to describe, but 
the proof that their complexes of end reductions have the stated form 
is rather lengthy. In the present paper we give a different method 
for constructing examples of \rirr\ Whitehead manifolds $W$ 
which cover \tm s $M$ with $\pi_1(M)\cong\mathbf{Z}$ and have \SW\ a 
triangulation of \R. This method has the advantage that the proof is much 
shorter. In addition we are able to classify all the end reductions of 
these examples. For the modified Scott-Tucker manifolds we were able 
to classify only those which are \rirr\ (which is sufficient to 
determine the complex). This gives the first \rirr\ Whitehead 
manifolds (other than those of genus one) for which the entire set of 
end reductions is known. 

The methods of this paper can also be used to construct \rirr\ 
Whitehead manifolds which cover \tm s with non-Abelian 
free fundamental groups and can cover only \tm s with free 
fundamental groups. This will be the subject of a later paper. 

The paper is organized as follows. In section 2 we give general background 
information and terminology. In section 3 we state those properties 
of end reductions we will need. In section 4 we prove the existence 
of graphs in the 3-ball having certain properties that we will need in our 
construction. In section 5 we prove the main technical result needed to 
determine the \er s of our examples. It is a condition on the 
embedding of one handlebody in the interior of another which ensures 
that any knot in the smaller handlebody which meets certain  
compressing disks for the boundary of the smaller handlebody in an essential 
way must meet all the compressing disks for the boundary of 
the larger handlebody. This result may be of some 
independent interest. In section 6 we give our basic construction of the 
examples $W$. In section 7 we prove some of their important properties. 
In section 8 we determine \SW. In section 9 we show how to modify 
the construction to get uncountably many such $W$.

%Section 2
\section{Background}

In general we follow \cite{He} or \cite{Ja} for basic \tm\ terminology. One 
slight difference is our use of the term \binc. This is usually 
reserved for surfaces $F$ which are properly embedded in a \tm\ 
$M$. We extend this to the case where $F$ is a surface in $\bd M$ as 
follows. $F$ is \textit{\binc} if whenever \De\ is a properly embedded 
disk in $M$ with $\De\cap F$ an arc \al\ and $\De\cap\overline{(\bd M-F)}$ 
an arc \be, then \al\ must be \bpl\ in $F$. 

When $X$ is a submanifold of $Y$ we denote the topological interior of 
$X$ by $\Inte X$ and the manifold interior of $X$ by $\inte X$. The 
\textit{exterior} of 
$X$ is the closure of the complement of a regular neighborhood of $X$ in 
$Y$. This term is also applied to the case of a graph $\Ga$ in $Y$. The 
regular neighborhood is denoted $N(\Ga,Y)$. A \textit{meridian} of an edge 
\ga\ of \Ga\ is 
the boundary of a properly embedded disk in $N(\Ga,Y)$ which meets \ga\ 
transversely in a single point.  

A sequence $\{C_n\}_{n\geq0}$ of compact, connected \tm s $C_n$ in a 
Whitehead manifold $W$ such that $C_n\sbs \inte C_{n+1}$ and 
$W-\inte C_n$ has no 
compact components is called a \textit{quasi-exhaustion in} $W$. 
If $\cup_{n\geq0} C_n=W$, then it is called an \textit{exhaustion for} $W$. 

The \textit{genus of} $\{C_n\}_{n\geq0}$ is the maximum of the genera of 
$\bd C_n$ or $\infty$ if these genera are unbounded. The \textit{genus of} 
$W$ is the minimum of the genera of its exhaustions. 

A plane $\Pi$ in $W$ is \textit{proper} if for each compact $K\sbs W$ 
one has that 
$K\cap\Pi$ is compact. A proper plane $\Pi$ is \textit{trivial} if some 
component of 
$W-\Pi$ has closure \hm\ to \halfspace. $W$ is \textit{\rirr} if every 
proper plane in $W$ is trivial. Every genus one Whitehead manifold 
is \rirr\ \cite{Kn}. 

A compact \tm\ $Y$ is \textit{weakly \anan} if every properly embedded 
\inc\ annulus in $Y$ has its boundary in a single component of $\bd Y$. 

%Lemma 2.1
\begin{lem} Suppose that for each compact $K\sbs W$ there is a 
quasi-exhaustion for $W$ such that 
\begin{enumerate} 
\item each $C_n$ is \irr,  
\item each $\bd C_n$ is \inc\ in $W-\inte C_n$, 
\item each $C_{n+1}-\inte C_n$ is \irr, \birr, and weakly \anan, and 
\item  $K\sbs C_1$. 
\end{enumerate}
Then $W$ is \rirr. \end{lem}

\begin{proof} This is Lemma 10.3 of \cite{My Endcov}, which derives from 
Lemma 4.2 of Scott and Tucker \cite{ST}. \end{proof}

%Section 3
\section{End reductions} 

In this section we collect some information about end reductions and 
define the complex of end reductions \SW\ of a Whitehead manifold $W$. 

A compact, connected 3-manifold $J$ in $ W$ is 
\textit{regular in} $W$ if $W-J$ is irreducible and has no component with 
compact closure. Since $W$ is \irr\ the first condition is equivalent to the 
statement that $J$ does not lie in a 3-ball in $W$. A quasi-exhaustion 
$\{C_n\}_{n\geq0}$ in $W$ is \textit{regular} if each $C_n$ is regular in $W$. 

Let $J$ be a regular 3-manifold in $W$, and let $V$ be an open subset 
of $W$ which contains $J$. We say that $V$ is \textit{\ei\ rel $J$ in 
$W$} if there is a regular \qe\ $\{C_n\}_{n\geq0}$ in $W$ such that 
$V=\cup_{n\geq0}C_n$, $J=C_0$, and $\bd C_n$ is \inc\ in $W-\inte J$ for all 
$n\geq0$. We say that $V$ has the \textit{engulfing property rel $J$ in $W$} 
if whenever $N$ is regular in $W$, $J\sbs\inte N$, and $\bd N$ is \inc\ in 
$W-J$, then $V$ is ambient isotopic rel $J$ to $V\p$ such that $N\sbs V\p$. 
$V$ is an \textit{\er\ of $W$ at $J$} if $V$ is \ei\ rel $J$ in $W$, $V$ 
has the engulfing property rel $J$ in $W$, and no component of $W-V$ has 
compact closure. 

%Theorem 3.1
\begin{thm}[Brin-Thickstun] Given  a regular 3-manifold $J$ in $W$, 
an \er\ $V$ of $W$ at $J$ exists and is unique up to non-ambient 
isotopy rel $J$ in $W$. \end{thm}

\begin{proof} This constitutes Theorems 2.1 and 2.3 of \cite{BT}. \end{proof}

It may help the reader's intuition about $V$ to see a brief sketch of its 
construction. We begin with a regular exhaustion $\{K_n\}_{n\geq0}$ of $W$ 
with $K_0=J$. Set $K_0^*=K_0$. If $\bd K_1$ is \inc\ in $W-J$ set $K_1^*=K_1$. 
Otherwise we ``completely compress'' $\bd K_1$ in $W-K_0^*$ to obtain $K_1^*$. 
We may assume that $K_1^*\sbs\inte K_2$. If $\bd K_2$ is \inc\ in $W-J$ we 
set $K_2^*=K_2$. Otherwise we completely compress $\bd K_2$ in $W-K_1^*$ to 
get $K_2^*$. We continue in this fashion to construct a sequence 
$\{K_n^*\}_{n\geq0}$. We let $V^*=\cup_{n\geq0}K_n^*$ and then let $V$ be the 
component of $V^*$ containing $J$. 

%Proposition 3.2
\begin{prop} Let $V$ be an \er\ of $W$ at $J$. Then the following hold: 
\begin{enumerate}
\item (Brin-Thickstun) If $J\p$ is regular in $W$, $J\sbs\inte J\p$, 
$J\p\sbs V$, and $\bd J\p$ is \inc\ in $W-J$, then $V$ is an \er\ of $W$ at 
$J\p$.  
\item There is a knot $\ka$ in $\inte J$ such that $V$ is an end reduction 
of $W$ at (a regular neighborhood of) \ka. 
\item $V$ is a Whitehead manifold. 
\end{enumerate}
\end{prop}

\begin{proof} (1) is Corollary 2.2.1 of \cite{BT}. (2) is Lemma 2.4 of 
\cite{My Endcov}. (3) is Lemma 2.6 of \cite{My Endcov}. \end{proof}

An \er\ $V$ of $W$ at $J$ is \textit{minimal} if whenever $U$ is an end 
reduction of $W$ at $K$ and $U\sbs V$, then there is a non-ambient isotopy 
of $U$ to $V$ in $W$. It is easily seen that genus one end reductions are 
minimal; recall that they are also \rirr. 

In \cite{Tu} Tucker constructed a \tm\ $W_0$ whose interior and boundary 
are \hm, respectively, to \RRR\ and \RR\ but which is not \hm\ to \halfspace. 
$W_0$ is a monotone union of solid tori which meet $\bd W_0$ in a monotone 
union of disks. It can be shown that the double of $W_0$ along its boundary 
is a Whitehead manifold which is a minimal end reduction of itself but is 
not \rirr. 

In \cite{My Endcov} and this paper examples are given of \rirr\ Whitehead 
manifolds having \rirr\ \er s which are not minimal. 

If $V$ is an \er\ of $W$, then we denote the non-ambient 
isotopy class of $V$ in $W$ by $[V]$. These isotopies are not required to 
be rel $J$. From now on we will usually drop the phrase ``non-ambient'' 
from ``non-ambient isotopy''. The vertices of \SW\ are those $[V]$ for 
which $V$ is minimal and \rirr. 

Distinct vertices $[V_0]$ and $[V_1]$ are joined by an edge if there is an 
\rirr\ \er\ $E_{0,1}$ of $W$ such that (1) $E_{0,1}$ contains 
representatives of $[V_0]$ and $[V_1]$, (2) every \rirr\ \er\ of $W$ 
contained in $E_{0,1}$ is isotopic in $W$ to $V_0$, $V_1$, or $E_{0,1}$, and 
(3) $[E_{0,1}]$ is unique among \rirr\ \er s of $W$ with respect to (1) and 
(2). 

Three distinct vertices $[V_0]$, $[V_1]$, and $[V_2]$ span a 2-simplex of 
\SW\ if each pair of vertices is joined by an edge and there is an \rirr\ 
\er\ $T_{0,1,2}$ of $W$ such that (1) $T_{0,1,2}$ contains representatives 
of each $[V_i]$ and $[E_{i,j}]$, (2) every \rirr\ \er\ of $W$ contained in 
$T_{0,1,2}$ is isotopic in $W$ to one of the $V_i$ or $E_{i,j}$ or to 
$T_{0,1,2}$, (3) $[T_{0,1,2}]$ is unique among \rirr\ \er s of $W$ with 
respect to (1) and (2). 

There is an obvious generalization of these definitions which inductively 
defines simplices of higher dimensions. 

Let $Homeo(W)$ denote the group of self-homeomorphisms of $W$. Let 
$Aut(\SW)$ denote the group of simplicial automorphisms of \SW. 
Each $g\in Homeo(W)$ induces a $\ga\in Aut(\SW)$. Let $\Psi:Homeo(W)
\rightarrow Aut(\SW)$ be the homomorphism given by $\Psi(g)=\ga$.   

%Theorem 3.3
\begin{thm} If $W$ is a non-trivial covering space of a \tm\ 
$M$ with group of covering translations $G\cong\pi_1(M)$, then the 
restriction $\Psi|G:G\rightarrow Aut(\SW)$ is one to one. \end{thm}

\begin{proof} This is Theorem 17.1 of \cite{My Endcov}. \end{proof}

%Corollary 3.4
\begin{cor} If \SW\ is isomorphic to a triangulation of \R, then 
$\pi_1(M)\cong\mathbf{Z}$. \end{cor}

\begin{proof} $\pi_1(M)$ must be torsion-free. The only non-trivial 
torsion-free subgroups of the infinite dihedral group $Aut(\SW)$ are 
infinite cyclic. \end{proof}

%Section 4
\section{Some poly-excellent graphs in the 3-ball}

A compact, connected, orientable \tm\ is \textit{superb} if it is \irr, 
\birr, and \anan, it contains a two-sided, properly 
embedded \inc\ surface, and it is not a 3-ball. It is \textit{excellent} if, 
in addition, it is atoroidal. In this paper superb \tm s which are 
not excellent will occur only in the last section. A compact, properly 
embedded 1-manifold in a compact, connected, orientable \tm\ is 
\textit{superb} or \textit{excellent} if its exterior is, respectively, 
superb or excellent. 
It is \textit{poly-superb} or \textit{poly-excellent} if for each non-empty 
collection of  its components the union of that collection is, respectively, 
superb or excellent. 

Define a \textit{k-tangle} to be a disjoint union of $k$ properly embedded 
arcs in a 3-ball. 

%Lemma 4.1
\begin{lem} For all $k\geq1$ poly-excellent $k$-tangles exist. \end{lem}

\begin{proof} This is Theorem 6.3 of \cite{My Attach}. \end{proof}

In this section we generalize this to certain graphs in the 3-ball. 
For $n\geq2$ define an \textit{$n$-frame} $F$ to be a graph having one vertex 
of degree $n$ and $n$ vertices of degree one; thus it is the cone on a 
set of $n$ points. A \textit{subframe} of $F$ is a subgraph of $F$ which is an 
$m$-frame for some $m\geq2$. Note that a single edge of $F$ is not a 
subframe of $F$. 

$F$ is \textit{properly embedded} in a 3-ball $B$ if $F\cap\bd B$ is the set 
of vertices of $F$ of degree one. A \textit{system of frames in} $B$ is a 
disjoint union $\mathcal{F}$ of finitely many properly embedded $n_i$-frames 
$F_i$ in $B$. We say that \FF\ is \textit{superb} or \textit{excellent} if 
its exterior is, respectively, superb or excellent. It is 
\textit{\ps} or \textit{\pe} if 
every non-empty subgraph of \FF\ whose components are subframes of the 
components of \FF\ is, respectively, superb or excellent. 
Note that the subgraph need not meet every component of \FF. 

%Theorem 4.2
\begin{thm} Let $k\geq1$. Suppose $n_1\geq2$. If $k\geq2$ assume that 
$n_i=2$ for $2\leq i\leq k$. Then there exists a \pe\ system \FF\ of 
$n_i$-frames $F_i$ in the 3-ball $B$. \end{thm}

In this paper we will need only the case $n_1=3$, but it is no harder 
to prove for $n_1>3$. 

We will need the following lemma for gluing together superb or 
excellent \tm s to obtain a superb or, respectively, excellent \tm. 

%Lemma 4.3
\begin{lem} Let $Y$ be a compact, connected, orientable \tm. 
Let $S$ be a compact, properly embedded, two-sided surface in $Y$. Let $Y\p$ 
be the \tm\ obtained by splitting $Y$ along $S$. Let $S\p$ and $S\pp$ be the 
two copies of $S$ which are identified to obtain $Y$. If each component of 
$Y\p$ is superb (respectively excellent), $S\p$, $S\pp$, and 
$(\bd Y\p)-\inte(S\p\cup S\pp)$ are \inc\ in $Y\p$, and each component of 
$S$ has negative Euler characteristic, then $Y$ is superb (respectively 
excellent). \end{lem}

\begin{proof} In the excellent case this is Lemma 2.1 of \cite{My Excel}. 
The superb case follows from the proof of that lemma. \end{proof}

\begin{proof}[Proof of Theorem 4.2] By Lemma 4.1 we may assume that 
$n_1\geq3$. 

We first prove the case $k=1$. Let $n=n_1$. Let $(\rho,\theta,\phi)$, 
$\rho\geq 0$, $0\leq\theta\leq2\pi$, $0\leq\phi\leq\pi$, be spherical 
coordinates in \RRR. We regard $B$ as the set $\rho\leq2$. Let $B\p$ be 
the set $\rho\leq1$. Let $\Sigma$ be the spherical shell $B-\inte B\p$. The 
$n$ halfplanes $\theta=0,2\pi/n,\ldots,2\pi(n-1)/n$ meet $\Sigma$ in disks 
$D_0,D_1,\ldots,D_{n-1}$ whose union cuts $\Sigma$ into 3-balls 
$B_0,B_1,\dots,B_{n-1}$, where $\bd B_j=D_j\cup D_{j+1}\cup E_j\cup E\p_j$ 
(subscripts taken mod $n$), where $E_j=B_j\cap \bd B$ and 
$E_j\p=B_j\cap\bd B\p$. We may think of 
$\Sigma$ as a cantaloupe which has been cut into $n$ wedge shaped slices and 
whose seeds have been removed. See Figure 1 for a schematic diagram of 
the following construction. 

In each $B_j$ we choose a \pe\ $(n+1)$-tangle 
$\al_{j,0}\cup\al_{j,1}\cup\cdots\cup\al_{j,n}$. We require (taking the 
subscript $j$ mod $n$) that $\al_{j,0}$ runs 
from $\inte E_j$ to $\inte D_{j+1}$, $\al_{j,p}$ runs from $\inte D_j$ to 
$\inte D_{j+1}$ for $1\leq p\leq n-1$, and $\al_{j,n}$ runs from $\inte D_j$ 
to $\inte E_j\p$. In addition we require that 
$\al_{j,p}\cap D_{j+1}=\al_{j+1,p+1}\cap D_{j+1}$. We then let 
$\be_j=\al_{j,0}\cup \al_{j+1,1}\cup\cdots\cup\al_{j-1,n-1}\cup\al_{j,n}$. 
The $\be_j$ are disjoint arcs each of which joins $\bd B$ to $\bd B\p$ in 
$\Sigma$. We may think of regular neighborhoods of the $\be_j$ as tunnels 
eaten out of the cantaloupe by $n$ worms who start on the outside and eat 
their way to the seed chamber in such a way that they each wind all the 
way around the cantaloupe, passing through every slice from one side to 
the other while coordinating their movements so that the union 
of the tunnels in each slice is \pe. 

\begin{figure}
\epsfig{file=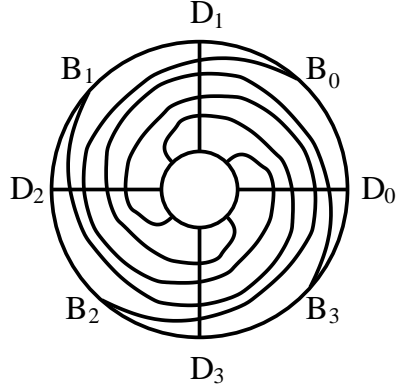, height=2in}
%\scalebox{0.5}{\includegraphics{fig1.pdf}}
\caption{The cantaloupe trick}
\end{figure}

The exterior in $\Sigma$ of the union of the $\be_j$ is equal to the 
exterior in $B$ of an $n$-frame $F$. We claim that $F$ is \pe. 
Let $F\p$ be an $m$-frame 
which is a subframe of $F$. Let $X_j$ be the exterior of $F\p\cap B_j$, 
and let $S_j=X_j\cap D_j$. Each $X_j$ is excellent. Since $m\geq2$ 
we have that $F\p$ meets 
each $D_j$ at least twice. Thus $\chi(S_j)<0$. Since no arc $\al_{j,p}$ 
joins $D_j$ to 
itself or $D_{j+1}$ to itself we have that $S_j$ and $S_{j+1}$ are each \inc\ 
in $X_j$. Since $X_j$ is \birr\ and neither $S_j$ nor $S_{j+1}$ is a disk we 
have that $\bd X_j-\inte S_j$ and $\bd X_j-\inte S_{j+1}$ are \inc\ in $S_j$. 
By successive applications of Lemma 4.3 we get that 
$X_0\p=X_1\cup\cdots\cup X_{n-1}$ is excellent. 
Now $X_0$ and $X_0\p$ are glued along the surface $S_0\cup S_1$, which 
is a disk with $2m+1$ holes. 
We may assume that $F\p$ meets $E_0$ and $E_0\p$. 
$\bd X_0-\inte(S_0\cup S_1)$ is the disjoint 
union of $m+1$ annuli. Since $X_0$ is \birr\ it follows that $S_0\cup S_1$ and 
$\bd X_0-\inte(S_0\cup S_1)$ are \inc\ in $X_0$. Now 
$\bd X_0\p-\inte(S_0\cup S_1)$ is the disjoint union of 
an annulus and two disks with $m-1$ holes. 
Since $X_0\p$ is \birr\ it follows that $S_0\cup S_1$ and 
$\bd X_0\p-\inte(S_0\cup S_1)$ are \inc\ in $X_0\p$. 
So by Lemma 4.3 the exterior $X_0\cup X_0\p$ of $F\p$ is excellent. 

We next prove the case $k>1$. We modify the construction of the 
previous case as follows. In $B_0$ we choose a \pe\ $(n+k)$-tangle 
$\al_{j,0}\cup\al_{j,1}\cup\cdots\cup\al_{j,n}\cup\ga_2\cup\cdots\cup\ga_k$, 
where each $\ga_q$ runs from $\inte E_0$ to 
itself. The $\al_{0,p}$ have the same properties as before. There is no change 
in the $B_j$ for $j\neq0$. Each $\ga_q$ is an arc and hence can be regarded 
as a 2-frame. The proof of poly-excellence works much as before. The only 
notable difference is that if the $n_1$-frame is deleted, then $B$ is the 
union of $B_0$ and a 3-ball along the disk $S_0\cup S_1\cup E_0\p$, and 
so $\ga_2\cup\cdots\cup\ga_k$ is 
\pe\ in $B$. \end{proof}

%Section 5
\section{Disk busting knots in handlebodies}

In this section we consider a knot \ka\ in the interior of a handlebody 
$C$ which is embedded in the interior of a handlebody \whC. We assume 
that $C$ and \whC\ each have genus at least one. Let \DD\ be a disjoint union 
of finitely many properly embedded disks in $C$ such that \DD\ splits $C$ 
into a collection of 3-balls and no component of \DD\ is \bpl\ in $C$. 
We say that \ka\ is \textit{\DD-busting} if no compressing disk for $\bd C$ 
in $C-\ka$ has 
the same boundary as a component of \DD. We give conditions on the 
embedding of $C$ in \whC\ which insure that if \ka\ is \DD-busting in $C$, 
then it 
is \textit{disk busting in} \whC, by which we mean that 
$\bd \whC$ is \inc\ in $\whC-\ka$. 

An \textit{$n$-pod} is a pair $(B,P)$ consisting of a 3-ball $B$ and a 
disjoint 
union $P$ of $n$ disks in $\bd B$. The components of $P$ are called the 
\textit{feet} 
of the $n$-pod. For $n=2$ or $n=3$ we use the term \textit{bipod} or 
\textit{tripod}, respectively. 

Two compact, properly embedded surfaces $S$ and $T$ in a \tm\ 
are in \textit{minimal general position} if they are in general position and 
among all such surfaces $S\p$ isotopic to $S$ one has that $S\cap T$ has the 
fewest components.   

%Lemma 5.1
\begin{lem} Let \whC\ be a handlebody of genus at least one. Let \whEE\ be 
a disjoint union of properly embedded disks in \whC\ which splits \whC\ into 
a union $(\whBB,\whPP)$ of bipods and tripods. Let \ka\ be a knot in 
$\inte \whC$ which is 
in general position with respect to \whEE. Let $(\ka\p,\bd\ka\p)$ be the 
1-manifold in 
$(\whBB,\whPP)$ obtained by splitting \ka\ along $\ka\cap\whEE$. Suppose 
that 
\begin{enumerate}
\item $\whPP-\ka\p$ is \inc\ in $\whBB-\ka\p$, 
\item $\whPP-\ka\p$ is \binc\ in $\whBB-\ka\p$, and 
\item each foot of $(\whBB,\whPP)$ meets $\ka\p$. 
\end{enumerate} 
Then \ka\ is disk busting in \whC. \end{lem}

\begin{proof} Suppose $D$ is a compressing disk for $\bd\whC$ in $\whC-\ka$. 
Put $D$ in \mgp\ with respect to $\whEE-\ka$. 

Suppose $D\cap(\whEE-\ka)$ contains a simple closed curve \ga. We may 
assume that \ga\ is innermost on $D$, so $\ga=\bd\De$ for a disk \De\ in $D$ 
with $\De\cap(\whEE-\ka)=\ga$. By (1) $\ga=\bd\Dep$ for a disk \Dep\ in 
$\whPP-\ka\p$. Then $\De\cup\Dep$ 
is a 2-sphere which bounds a 3-ball in \BB\ which by (3) misses $\ka\p$. Thus 
there is an isotopy of $D$ in $\whC-\ka$ which removes at least \ga\ from the 
intersection, thereby contradicting minimality. 

Now suppose that $D\cap(\whEE-\ka)$ has a component \al\ which is an arc. We 
may assume that \al\ is outermost on $D$, so there is an arc \be\ in $\bd D$ 
such that $\bd\al=\bd\be$ and $\al\cup\be=\bd\De$ for a disk \De\ in $D$ 
with $\De\cap(\whEE-\ka)=\al$. 
By (2) there is a disk \Dep\ in $\whPP-\ka\p$ and an arc $\al\p$ in $\bd\whPP$ 
such that 
$\al\cap\al\p=\bd\al=\bd\al\p$ and $\bd\Dep=\al\cup\al\p$. Then $\De\cup\Dep$ 
is a disk with 
$\bd(\De\cup\Dep)=\al\p\cup\be$. By (1) and (3) $\al\p\cup\be=\bd\Depp$ for 
a disk \Depp\ in $\bd\whBB-\Inte\whPP$. 
We have that $\De\cup\Dep\cup\Depp$ is a 2-sphere bounding a 3-ball in \whBB\ 
which 
by(3) misses $\ka\p$. Thus there is an isotopy of $D$ in $\whC-\ka$ which 
removes 
at least \al\ from the intersection, thereby contradicting minimality. 

We now have that $D\cap(\whEE-\ka)=\ns$, so $D$ lies in some component of 
$\whBB$. If $\bd D$ does not bound a disk in $\bd\whBB-\Inte\whPP$, then it 
is \pl\ in this 
surface to a component of $\bd\whPP$, thereby contradicting (1) and (3). 
\end{proof}
 
An $n$-pod $(B,P)$ is \textit{properly embedded} in an $m$-pod $(\whB,\whP)$ 
if $B\sbs\whB$ and $B\cap\bd\whB=B\cap\inte \whP=P$. Note that $(B,P)$ is a 
regular neighborhood of an $n$-frame in \whB. 

%Lemma 5.2
\begin{lem} Let $(\whB,\whP)$ be a bipod or tripod. Let $(\BB,\PP)$ be a 
disjoint union of bipods and tripods properly embedded in $(\whB,\whP)$. 
Let \la\ be a disjoint union of finitely many arcs properly embedded 
in \BB\ with $\la\cap\bd\BB=\bd\la\sbs\inte\PP$. Suppose that 
\renewcommand{\labelenumi}{(\roman{enumi})}
\begin{enumerate}
\item $\PP-\la$ is \inc\ in $\BB-\la$, 
\item $\PP-\la$ is \binc\ in $\BB-\la$, 
\item each foot of \PP\ meets \la\, 
\item each foot of \whP\ meets \PP, 
\item $\bd\BB-\inte\PP$ and $\whP-\inte\PP$ are \inc\ in $\whB-\Inte\BB$, and 
\item if $(\whB,\whP)$ is a tripod, then $\whP-\inte\PP$ is 
\binc\ in $\whB-\Inte\BB$, 
\item if any component of $(\BB,\PP)$ is a tripod, then $\whB-\Inte\BB$ is 
\birr. 
\end{enumerate}
Then 
\renewcommand{\labelenumi}{(\arabic{enumi})}
\begin{enumerate}
\item $\whP-\la$ is \inc\ in $\whB-\la$, 
\item $\whP-\la$ is \binc\ in $\whB-\la$, and 
\item each foot of \whP\ meets \la. 
\end{enumerate}
\end{lem}

\begin{proof} Suppose $D$ is a compressing disk for $\whP-\la$ in $\whB-\la$. 
Put $D$ in \mgp\ with respect to $\bd\BB-\inte\PP$. 

Suppose $D\cap(\bd\BB-\inte\PP)$ has a simple closed curve component \ga. 
We may assume that \ga\ is innermost on $D$, so $\ga=\bd\De$ for a disk \De\ 
in $D$ with $\De\cap(\bd\BB-\inte\PP)=\ga$. 

If \De\ lies in $\BB-\la$, then it follows from (i) and (iii) that 
$\ga=\bd\Dep$ for a disk \Dep\ in $\bd\BB-\inte\PP$. Then $\De\cup\Dep$ is a 
2-sphere which bounds a 3-ball in \BB\ which misses \la, so there is an 
isotopy of $D$ in $\whB-\la$ which removes at least \ga\ from the 
intersection, contradicting minimality. 

If \De\ lies in $\whB-\Inte\BB$, then by (v) there is a disk \Dep\ in 
$\bd\whB-\inte\PP$ such that $\ga=\bd\Dep$. Then $\De\cup\Dep$ is a 2-sphere 
which bounds a 3-ball in $\whB-\Inte\BB$ which misses \la, so there is an 
isotopy of $D$ in $\whB-\la$ which removes at least \ga\ from the 
intersection, contradicting minimality. 

Thus there are no simple closed curve components. Suppose there is 
a component \al\ which is an arc. We may assume that \al\ is outermost 
on $D$, so there is an arc \be\ in $\bd D$ such that $\bd\al=\bd\be$ and 
$\al\cup\be=\bd\De$ for a disk \De\ in $D$ with 
$\De\cap(\bd\BB-\inte\PP)=\al$. 

If \De\ lies in $\BB-\la$, then \be\ lies in $\PP-\la$. By (ii) there is a 
disk \Dep\ in $\PP-\la$ and an arc $\be\p$ in $\bd\PP$ such that 
$\be\cap\be\p=\bd\be=\bd\be\p$ and $\bd\Dep=\be\cup\be\p$. 
Then $\De\cup\Dep$ is a disk with $\bd(\De\cup\Dep)=\al\cup\be\p$. 
By (i) and (iii) there is a disk \Depp\ in $\bd\BB-\inte\PP$ with 
$\bd\Depp=\al\cup\be\p$. Then $\De\cup\Dep\cup\Depp$ is a 2-sphere 
which bounds a 3-ball in \BB\ that misses \la. Thus there is an isotopy of $D$ 
in $\whB-\la$ which removes at least \ga\ from the intersection, contradicting 
minimality. 

If \De\ lies in $\whB-\Inte\BB$, the \be\ lies in $\whP-\PP$. 

Suppose the component $(B,P)$ of $(\BB,\PP)$ containing \al\ is a bipod. 

Assume that \al\ joins the two feet of $(B,P)$. Since \be\ lies in $\whP-\PP$ 
these two feet must lie in the same foot of $(\whB,\whP)$. Let $N$ be a 
regular neighborhood of $B\cup\De$ in \whB. Then $N$ is a 3-ball such that 
$N\cap\bd\whB=N\cap\whP$ and is a disk $\whE$. The disk $\bd N-\inte\whE$ 
is therefore a compressing disk for $\whP-\inte\PP$ in $\whB-\Inte\PP$, 
contradicting (v). 

Thus $\bd\al$ lies in a single foot of $(B,P)$. 
Then there is a disk \Dep\ in $\bd B-\inte P$ with 
$\bd\Dep=\al\cup\al\p$, where $\al\p$ is an arc in $\bd P$ 
with $\bd\al=\bd\al\p$. 
So $\De\cup\Dep$ is a disk in $\whB-\Inte\BB$ with 
$\bd(\De\cup\Dep)=\al\p\cup\be$. By (v) there is a disk \Depp\ in 
$\whP-\inte\PP$ with $\bd\Depp=\al\p\cup\be$. Then $\De\cup\Dep\cup\Depp$ is 
a 2-sphere bounding a 3-ball in \whB\ which misses \la. Thus there is an 
isotopy of $D$ in $\whB-\ka$ which removes at least \al\ from the 
intersection, contradicting minimality. 

Suppose the component $(B,P)$ of $(\BB,\PP)$ containing \al\ is a tripod. 
By (vii) there is a disk \Dep\ in $\bd(\whB-\Inte\BB)$ such that 
$\bd\Dep=\bd\De$. Since each component of $\bd\PP$ is a non-separating 
curve on $\bd(\whB-\Inte\BB)$ we must have that $\bd\be$ lies in a single 
component of $\bd\PP$. Moreover \Dep\ is the union of a disk in 
$\bd\BB-\Inte\PP$ and a disk in $\whP-\inte\PP$ which meet along an arc 
in $\bd\PP$, and $\De\cup\Dep$ is a 2-sphere bounding a 3-ball in \whB\ 
which misses \la. Thus there is an isotopy of $D$ in $\whB-\la$ which removes 
at least \al\ from the intersection, contradicting minimality. 

So we have that $D$ misses $\bd\BB-\inte\PP$. If $D$ lies in $\BB-\la$, 
then by (i) $\bd D=\bd D\p$ for a disk $D\p$ in $\PP-\la$. If $D$ lies in 
$\whB-\Inte\BB$, then by (v) $\bd D=\bd D\p$ for a disk $D\p$ in 
$\whP-\inte\PP$. This completes the proof of (1). 

Now suppose that $D$ is a $\bd$-compressing disk for $\whP-\la$ in $\whB-\la$. 
We have that $\bd D=\ga\cup\de$ for arcs \ga\ in $\whP-\la$ and \de\ in 
$\bd\whB-\inte\whP$. Put $D$ in \mgp\ with respect to $\bd\BB-\inte\PP$. 
As in the proof of (1) we may assume that no component of the intersection 
is a simple closed curve.

Suppose the intersection has a component \al\ which is an arc. We may 
assume that \al\ is outermost with respect to \de, by which we mean that 
there is a disk \De\ in $D$ and an arc \be\ in \ga\ wuch that 
$\bd\al=\bd\be$, $\bd\De=\al\cup\be$, and $\De\cap(\bd\BB-\inte\PP)=\al$. The 
analysis of \De\ now proceeds as in the proof of (1), and we again contradict 
minimality. 

So $D$ misses $\bd\BB-\inte P$, and $D$ lies in $\whB-\Inte\BB$. 
By (v) and (vi) $\bd D=\bd D\p$ for a disk $D\p$ in $\bd(\whB-\Inte\BB)$. 
Since each component of $\bd\PP$ is non-separating in $\bd(\whB-\Inte\BB)$ 
we have that $D\p\cap\PP$ is a disk. This completes the proof of (2). 

(3) follows from (iii) and (iv). \end{proof}

A disjoint union of $n_i$-pods $(B_i,P_i)$ properly embedded in an $m$-pod 
$(\whB,\whP)$ is \textit{\ps} or \textit{\pe} if the corresponding 
union of $n_i$-frames is, respectively, \ps\ or \pe. 

We suppose now that $\whEE$ is a disjoint union of properly embedded 
disks in \whC\ which splits \whC\ into a union $(\whBB,\whPP)$ of bipods 
and tripods. These bipods and tripods and their feet are called \textit{big}. 
We assume that $\whEE\cap C$ is a union \EE\ of properly embedded disks in 
$C$ which splits $C$ into a union $(\BB,\PP)$ of bipods and tripods. These 
bipods and tripods and their feet are called \textit{small}.  
Let $\DD$ be a union of components of $\EE$. Suppose \ka\ is a knot in 
$\inte C$ which is $\DD$-busting. 
A small foot is called \textit{hot} if it is 
\pl\ in $C$ to a component of \DD. It is \textit{warm} if there is no 
compressing disk for $\bd C$ in $C-\ka$ which has the same boundary. 
It is \textit{cold} if there is such a compressing disk. Note that 
every hot foot is warm. 

%Proposition 5.3
\begin{prop} Suppose that for each big bipod or tripod $(\whB,\whP)$  
\begin{enumerate}
\item each big foot of $(\whB,\whP)$ contains a small warm foot of 
$(\whB,\whP)\cap(\BB,\PP)$, and 
\item either 
\begin{enumerate}
\item $(\whB,\whP)\cap(\BB,\PP)$ is \ps, or 
\item $(\whB,\whP)$ is a bipod, $(\whB,\whP)\cap(\BB,\PP)$ consists of bipods, 
and each of these small bipods meets each of the two big feet of 
$(\whB,\whP)$. 
\end{enumerate}
\end{enumerate}
Then every \DD-busting knot \ka\ in $C$ is disk busting in \whC. \end{prop}

\begin{proof} Suppose \ka\ is \DD-busting in $C$. Isotop \ka\ in $C$ so that 
it is in \mgp\ with respect to \EE. We will show that after possibly modifying 
$(\BB,\PP)$ we will have that for each big bipod or tripod $(\whB,\whP)$ 
it is the case that $(\whB,\whP)\cap(\BB,\PP)$ satisfies the hypotheses of 
Lemma 5.2 with $\la=\ka\cap\whB$. Note that we do not require that the 
components of the modified $(\BB,\PP)$ match up along \whEE\ to give a new 
handlebody in \whC. 

So let $(\whB,\whP)$ be a big bipod or tripod. 

Suppose we are in case 2(a). 

Consider a small bipod $(B,P)$ in $(\whB,\whP)\cap(\BB,\PP)$. If 
$\la\cap B=\ns$, then we discard $(B,P)$ from $(\BB,\PP)$ to obtain a new 
\ps\ system. If $\la\cap B\neq\ns$, then by minimality \la\ meets each 
small foot of $(B,P)$ and 
$P-\la$ is \inc\ in $B-\la$. Since $(B,P)$ is a bipod we then have 
that $P-\la$ is \binc\ in $B-\la$. 

Consider a small tripod $(B,P)$. If $\la\cap B=\ns$, then we discard 
$(B,P)$ from $(\BB,\PP)$ to obtain a new \ps\ system. If $\la\cap B\neq\ns$, 
then by minimality \la\ meets at least two small feet of $(B,P)$. 

Suppose \la\ misses the third small foot. Then we push that foot slightly into 
$\inte \whB$ to obtain a bipod. This gives a new \ps\ system. We have that 
\la\ meets each foot of the new $(B,P)$, and $P-\la$ is \inc\ and \binc\ in 
$B-\la$. 

Suppose \la\ meets the third small foot. Then $P-\la$ is \inc\ in $B-\la$. 
If $P-\la$ is $\bd$-compressible 
in $B-\la$, then there is a properly embedded disk \De\ 
in $B-\la$ which meets a component $E$ of $P$ in an arc \al\ and 
$\bd B-\inte P$ in an arc \be\ such that $\bd\al=\bd\be$, $\bd\De=\al\cup\be$, 
and \al\ splits $E$ into two disks each of which meets \la. Since $E-\la$ 
is \inc\ in $B-\la$ we must have that the two components of 
$P-E$ are separated from each other by \De. We split $(B,P)$ along 
$(\De,\al)$ to obtain two bipods $(B\p,P\p)$ and $(B\pp,P\pp)$. We have that 
$(P\p\cup P\pp)-\la$ is \inc\ and \binc\ in $(B\p\cup B\pp)-\la$. The 
exterior of the new system is \hm\ to that of the old system by a 
homeomorphism which is the identity on the other components of 
$\bd\BB-\inte\PP$, and so the new system is also \ps. 

The feet discarded by our modifications are precisely the cold feet of 
$(\whB,\whP)\cap(\BB,\PP)$. Some warm feet may be split into pairs of warm 
feet. It follows that conditions (i), (ii), and (iii) of Lemma 5.2 are 
satisfied. Since each component of \whP\ contains a warm foot condition (iv) 
is satisfied. Since our modifications preserve \ps ness condtions (v), (vi),  
and (vii) are also satisfied. 

Now suppose that we are in case 2(b). As in the previous case we 
discard all small bipods with cold feet and get that conditions (i), (ii), 
(iii), and (iv) are satisfied. Since each small bipod joins the two big 
feet condition (v) is satisfied. Conditions (vi) and (vii) are 
vacuously satisfied. 

The result now follows from Lemmas 5.1 and 5.2. \end{proof}

%Section 6
\section{The construction of $W$}

\begin{figure}
\epsfig{file=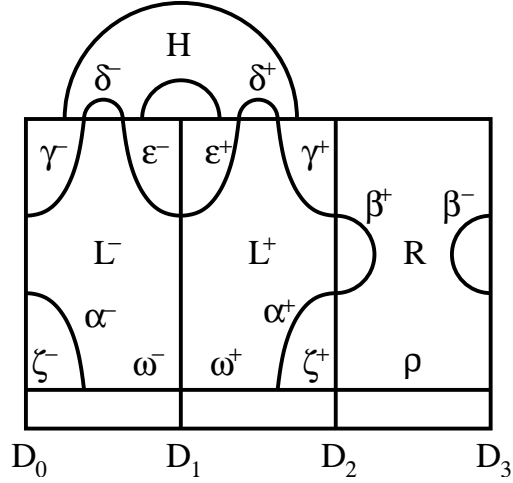, height=2.5in}
%\scalebox{0.5}{\includegraphics{fig2.pdf}}
\caption{The systems of frames in the pieces of $J$}
\end{figure}

In this section we construct an \rirr\ contractible open \tm\ $W$ 
which covers a \tm\ $W^{\#}$ with $\pi_1(W^{\#})\cong\mathbf{Z}$. 
It will be shown that \SW\ is a triangulation of \R\ and hence every 
\tm\ non-trivially covered by $W$ must have fundamental group $\mathbf{Z}$. 

Let $P=D\times[0,3]$, where $D$ is a closed disk. Let $L^-=D\times[0,1]$, 
$L^+=D\times[1,2]$, $R=D\times[2,3]$, and $D_j=D\times\{j\}$ for 
$j=0,1,2,3$. Let $L=L^-\cup L^+$. Attach a 1-handle $H$ to $P$ so that it 
joins $\bd D\times(0,1)$ to $\bd D\times(1,2)$, thus giving a solid torus 
$J=P\cup H$. Let $J^{\#}$ be the genus two handlebody obtained from $J$ by 
identifying $D_0$ and $D_3$. Let $P\sh$ be the solid torus in $J\sh$ which 
is the image of $P$ under the identification. With the exceptions of $J$, 
$J\sh$, $P$, and $P\sh$ we will usually use the same symbol for subsets of 
$J$ and their images in $J\sh$, relying on the context for the meaning. 
Thus we write, for example, $J\sh=P\sh\cup H$. 

We next define a certain graph \tht\ in $J\sh$ as follows. See Figure 2 
for a schematic diagram of this construction.

Choose a \ps\ system of frames in $L^-$ consisting of a 3-frame 
and two 2-frames. The 3-frame consists of arcs $\al^-$, $\zeta^-$, and 
$\omega^-$ meeting in a common endpoint in $\inte L^-$. The other endpoints 
of $\al^-$ and $\ze^-$ lie in $\inte D_0$. The other endpoint of $\om^-$ 
lies in $\inte D_1$. One 2-frame is an arc $\ga^-$ joining $\inte D_0$ and 
$\inte(L^-\cap H)$. The other 2-frame is an arc $\ep^-$ joining $\inte D_1$ 
and $\inte(L^-\cap H)$. 

Let $r$ be the homeomorphism $r(x,t)=(x,2-t)$ from $D\times[0,2]$ to itself 
which reflects in the disk $D_1$. We have that $r(L^-)=L^+$. Denote 
$r(\al^-)$ by $\al^+$, etc. This defines a \ps\ system of frames in $L^+$. 

Next choose a \ps\ 2-tangle in $H$ with components $\de^-$ and 
$\de^+$ such that $\bd\de^{\pm}=(\ga^{\pm}\cup \ep^{\pm})\cap H$. 
Then choose a \ps\ 3-tangle in $R$ with components $\be^-$, $\be^+$, and 
$\rho$, where $\bd\be^{\pm}=(\al^{\pm}\cup\ga^{\pm})\cap R$ and 
$\bd\rho=(\ze^-\cup\ze^+)\cap R$. 

Let $\eta$ be the arc 
$\al^-\cup\be^-\cup\ga^-\cup\de^-\cup\ep^-\cup\ep^+\cup\de^+\cup\ga^+\cup
\be^+\cup\al^+$, $\la$ the arc $\om^-\cup\om^+$, and $\mu$ the arc 
$\ze^+\cup\rho\cup\ze^-$. Set $\tht=\eta\cup\la\cup\mu$. 

For each integer $n\geq0$ take a copy of each of these objects. Denote 
the $n^{th}$ copy of $D_j$ by $D_{n,j}$, that of each of the other 
objects by a subscript $n$. We regard the arcs and graphs with subscripts 
$n$ as embedded in the \tm s with subscript $n+1$. 

We embed $J_n\sh$ in $\inte J_{n+1}\sh$ as follows. $L_n$ is sent to 
$N(\la_n,L_{n+1})$. $R_n$ is sent to $N(\mu_n\cap(P_{n+1}\sh-\inte L_n), 
P_{n+1}\sh-\inte L_n)$. $H_n$ is sent to $N(\eta_n\cap(J_{n+1}\sh-\inte P_n), 
J_{n+1}\sh-\inte P_n)$. 

Now let $W\sh$ be the direct limit of the $J_n\sh$, and let $p:W\rightarrow 
W\sh$ be the universal covering map. Then $\pi_1(W\sh)$ is infinite cyclic. 
Let $h:W\rightarrow W$ be a generator of the group of covering translations. 
We regard $p\n(P_n\sh)$ as $D_n\times\R$ with $P_{n,j}=D_n\times[3j,3j+3]$, 
$L_{n,j}^-=D_n\times[3j,3j+1]$, $L_{n,j}^+=D_n\times[3j+1,3j+2]$, and 
$R_{n,j}=D_n\times[3j+2,3j+3]$. We set $L_{n,j}=L_{n,j}^-\cup L_{n,j}^+$. 
We have that $p\n(H_n)$ is a disjoint union of 1-handles $H_{n,j}$, 
where $H_{n,j}$ is attached to $\bd D_n\times(3j,3j+2)$, thereby yielding 
a copy $J_{n,j}=P_{n,j}\cup H_{n,j}$ of $J_n$. Set $D_{n,k}=D_n\times\{k\}$ 
for $k\in\mathbf{Z}$. For all the objects with subscript $n$ contained in 
$J_{n+1}$ denote the component of the preimage contained in $J_{n+1,j}$ by 
the subscripts $n,j$. 
We denote by the same symbols $\eta_n$ and $\mu_n$ the components of the 
preimages of $\eta_n$ and $\mu_n$ which meet $\omega_n^+$. 
We assume that $h$ is chosen so that $h(D_{n,k})=D_{n,k+3}$ and 
the image under $h$ of any other object with subscripts $n,j$ has subscripts 
$n,j+1$. 

We next describe certain families of \qe s in $W$. Let $\PP=\{p_1,p_2,\ldots, 
p_m\}$ be a finite non-empty set of distinct integers with $p_1<p_2< \cdots
<p_m$. We say that \PP\ is \textit{good} if its elements are consecutive. 
Otherwise \PP\ is \textit{bad}. If $m=1$, then \PP\ is automatically good. 

For $n\geq0$ let $C_n^{\PP}$ be the union of those $R_{n,j}$ with 
$p_1-1\leq j\leq p_m$, those $L_{n,j}$ with $p_1\leq j\leq p_m$, and those 
$H_{n,p}$ with $p\in\PP$. Each $C_n^{\PP}$ is a cube with $m$ handles 
embedded in $\inte C_{n+1}^{\PP}$. In Figure 3 we give a schematic diagram 
for the case of $\PP=\{p,p+1,p+2\}$. 

The \qe\ $\{C_n^{\PP}\}_{n\geq0}$ is denoted by $C^{\PP}$; its union is 
denoted by $V^{\PP}$. Whenever \PP\ is good and $m>1$ we denote $V^{\PP}$ by 
$V^{p,q}$, where $p=p_1$ and $q=p_m$. When $\PP=\{p\}$ we use the notation 
$V^p$. The expressions $C_n^{p,q}$, $C_n^p$, $C^{p,q}$, and $C^p$ are 
defined similarly. 

\begin{figure}
\epsfig{file=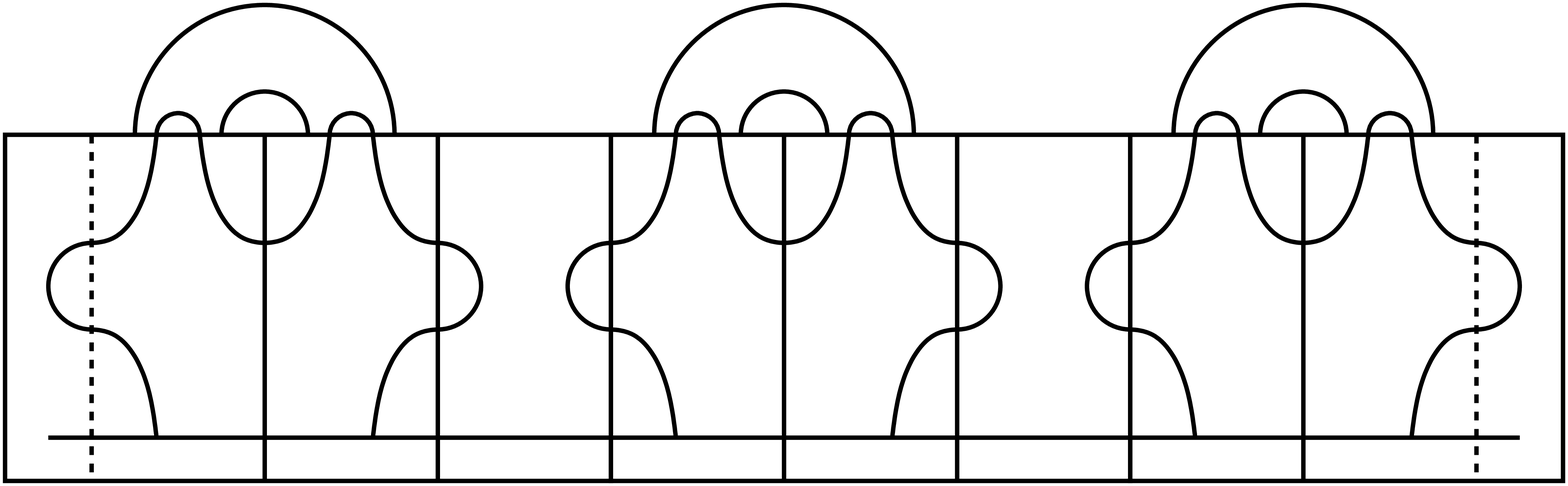, width=5in}
%\scalebox{0.25}{\includegraphics{fig3.pdf}}
\caption{The embedding of $C^{\PP}_n$ in $C^{\PP}_{n+1}$ for a good \PP}
\end{figure}

%Section 7
\section{Some properties of $W$}

Given $\PP=\{p_1,\ldots,p_m\}$ and $n>0$, let 
$Y=C^{\PP}_{n+1}-\inte C^{\PP}_n$, $p=p_1$, and $q=p_m$. If $m>1$ set 
$Z^-=Y\cap(R_{n+1,p-1}\cup L^-_{n+1,p})$, 
$Z^+=Y\cap(L^+_{n+1,q}\cup R_{n+1,q})$, 
$Z=Z^-\cup Z^+$, and $X=Y-\Inte Z$. 

%Lemma 7.1
\begin{lem} $Y$ is \irr\ and \birr. \end{lem}

\begin{proof} First consider the case $m=1$. Then $C^p_n$ is a solid torus 
in $C^p_{n+1}$ with winding number zero. Any compressing disk for 
$\bd C_{n+1,p}$ in $Y$ would be a meridinal disk for $C_{n+1,p}$. Since 
$\de^-_{n,p}\cup\de^+_{n,p}$ is \ps\ in $H_{n+1,p}$ we have that 
$H_{n,p}\cap L_{n,p}\cap Y$ is \inc\ in $H_{n+1,p}\cap Y$. 
It is \inc\ in $(C_{n+1,p}-\Inte H_{n+1,p})\cap Y$ for homological reasons. 
Thus $H_{n+1,p}\cap L_{n+1,p}\cap Y$ is \inc\ in $Y$ and thus so is 
$\bd C_{n+1,p}$. If $\bd C^p_n$ is compressible in $Y$, then the union 
of $C^p_n$ and a 2-handle with core the compressing disk is a 3-ball in 
$C^p_{n+1}$, and so $\bd C^p_{n+1}$ is compressible in $Y$, a contradiction. 

Now suppose $m>1$. Consider the surfaces $D_{n+1,k}\cap Y$ for 
$3p+1\leq k\leq 3q+1$ and $H_{n+1,p}\cap L^{\pm}_{n+1,p}\cap Y$ for $p\in\PP$. 
They split $Y$ into \irr\ pieces. With the exception of $Z^{\pm}$ it follows 
from \ps ness that each of these pieces is superb, and so each of those 
surfaces contained in its boundary is \inc\ and \binc. It follows that $X$ is 
\irr\ and \birr. $Z^{\pm}\cap X$ consists of two disks with two holes, and 
$\bd Z^{\pm}-\inte(Z^{\pm}\cap X)=\bd(Z^{\pm}\cap X)\times[0,1]$. 
Thus $Z^{\pm}\cap X$ is \inc\ and \binc\ in $Z^{\pm}$. Thus the result 
follows. \end{proof}

%Lemma 7.2
\begin{lem}$\bd C^{\PP}_n$ is \inc\ in $W-\inte C^{\PP}_n$. \end{lem}

\begin{proof} A compressing disk $D$ must lie in $C^{r,s}_m-\inte C^{\PP}_n$ 
for some $r\leq p$, $s\geq q$, and $m>n$. We can isotop $D$ off compressing 
disks for $\bd C^{r,s}_m$ in $C^{r,s}_m-\inte C^{\PP}_n$ so that it lies in 
$C^{\PP}_m-\inte C^{\PP}_n$. The result then follows from the previous 
lemma. \end{proof}

%Lemma 7.3
\begin{lem} If \PP\ is good and $m>1$, then 
\begin{enumerate}
\item if $A$ is an \inc\ annulus in $Y$, then $\bd A=\bd A\p$ for an 
annulus $A\p$ in $\bd Y$, and 
\item if $T$ is an \inc\ torus in $Y$, then $T$ bounds a compact \tm\ in $Y$. 
\end{enumerate}
\end{lem}

\begin{proof} (1) Put $A$ in \mgp\ with respect to $X\cap Z$. 
Let \al\ be a component of $A\cap X\cap Z$. Then \al\ is not 
a simple closed curve bounding a disk in $A$. 

Suppose \al\ is an outermost arc on $A$, so $\bd\De=\al\cup\be$ for an 
arc \be\ in $\bd A$ and a disk \De\ in $A$ with $\De\cap X\cap Z=\al$. 
If $\De\sbs X$, then $\bd\De=\bd\Dep$ for a disk $\Dep$ in $\bd X$. 
then $\De\cup\Dep$ bounds a 3-ball in $X$, and an isotopy across it removes 
at least \al\ from the intersection. If $\De\sbs Z$, then \be\ is 
\bpl\ in one of the annuli comprising $\bd Z-\inte(X\cap Z)$; it follows 
that one can again reduce the intersection. Thus \al\ is not an arc. 

So \al\ is a simple closed curve. $\bd A\p=\al\cup\be$ for some annulus 
component $A\p$ of $A\cap X$ and some \be\ in $(A\cap X\cap Z)\cup \bd A$. 
Then $A\p$ is \pl\ in $X$ to an annulus $A\pp$ in $\bd X$. If $A\pp$ lies in 
$X\cap Z$, then we can isotop to remove at least $\al\cup\be$. If $A\pp$ does 
not lie in $X\cap Z$, then either we can isotop to remove \al\ or $A\pp$ 
contains an annulus component $G$ of $\bd X-\inte(X\cap Z)$. We may assume 
that the centerline of $G$ is a meridian of $\be^+_{n,q}$ and that the 
component of $X\cap Z$ containing $\bd A\p$ is $F=H_{n+1,q}
\cap L^+_{n+1,q}\cap Y$. We may further assume that all the components of 
$A\cap X$ are \pl\ to $G$ and lie in $H_{n+1,q}\cap Y$. For homological 
reasons all the components of $A\cap Z$ must have their boundaries in the 
union of $F$ and the two annulus components of $\bd C^{\PP}_n\cap Z^+$. 
In particular, $\bd A$ lies in the union of these two annuli and so bounds an 
annulus in their union with $G$. 

Suppose $A\cap X\cap Z=\ns$. If $A\sbs X$, then $A$ is \pl\ in $X$ to an 
annulus $A\p$ in $\bd X$ with $\bd A\p$ in $\bd X-\inte(X\cap Z)$. It follows 
that $A\p$ lies in $\bd Y$. If $A\sbs Z$, then for homological reasons $\bd A$ 
must lie in one of the three annulus components of 
$\bd Z^{\pm}-\inte(X\cap Z^{\pm})$. 

(2) Suppose $T$ is in \mgp\ with respect to $X\cap Z$. 
$T$ cannot lie in $X$ since it would be \bpl\ in $X$, but $\bd X$ has no 
tori. If $T$ lies in $Z^{\pm}$, then since $\bd Z^{\pm}$ is connected $T$ 
must bound a compact \tm\ in $Z^{\pm}$. 

So we may assume that $T\cap Z\neq\ns$. Let $A$ be a component of $T\cap X$. 
As in the proof of (1) we may assume that $A$ is \pl\ in $X$ to an annulus 
$A\p$ in $\bd X$ which contains an annulus component $G$ of 
$\bd X-\inte(X\cap Z)$ whose centerline is a meridian of $\be^+_{n,q}$ and 
that all such components are \pl\ to $G$ and lie in $H_{n+1,q}\cap Y$. 
All the components of $T\cap Z^{\pm}$ must have their boundaries in the 
component $F$ of $X\cap Z^{\pm}$ which meets $G$. So $T$ lies in 
$(H_{n+1,q}\cup L^+_{n+1,q}\cup R_{n+1,q})\cap Y$. Since this \tm\ has 
connected boundary $T$ must bound a compact \tm\ in its interior. \end{proof}

%Lemma 7.4
\begin{lem} $V^p$ does not embed in \RRR. \end{lem}

\begin{proof} Since $\be^+\cup\be^-\cup\rho$ is \ps\ in $R$ we have that 
$\be^+_{n,p}$ is knotted in $R_{n,p}$. The result then follows from 
\cite{Ha}. \end{proof}

%Proposition 7.5
\begin{prop} $W$ is \rirr. If \PP\ is good, then $V^{\PP}$ is \rirr. 
\end{prop}

\begin{proof} It suffices to show that for each good \PP\ the \qe\ $C^{\PP}$ 
of $W$ satisfies conditions (1)--(3) of Lemma 2.1. When $m=1$ this 
follows from \cite{Kn}, so assume $m>1$. Each $C^{\PP}_n$ is a cube 
with handles, so is \irr. We have that $\bd C^{\PP}_n$ is \inc\ in 
$W-\inte C^{\PP}_n$ and that $Y$ is \birr\ and weakly \anan. \end{proof}

%Proposition 7.6
\begin{prop} If \PP\ is bad, then $V^{\PP}$ is not \rirr. \end{prop}

\begin{proof} There is an $s$ such that $p<s<q$ and $s\notin\PP$. We may 
assume that the embedding of $J\sh_n$ in $J\sh_{n+1}$ is such that 
$D_{n,1}\sbs D_{n+1,1}$ for all $n\geq0$. Then 
$D_{n,3s+1}\sbs\inte D_{n+1,3s+1}$ for all $n\geq0$. The union $\Pi$ of 
these disks is a plane which is proper in $V^{\PP}$ (but not in $W$!). 
$V^{\PP}-\Pi$ has two components, one containing $V^p$ and the other 
containing $V^q$. Since $V^p$ and $V^q$ do not embed in \RRR\ we have that 
$\Pi$ is non-trivial in $V^{\PP}$. \end{proof}

A \textit{classical knot space} is a space \hm\ to the exterior of a 
non-trivial knot in $S^3$. 

%Lemma 7.7
\begin{lem} If \PP\ is good and $m>1$, then every \inc\ torus $T$ in 
$V^{\PP}-\inte C^{\PP}_n$ bounds a compact \tm\ in $V^{\PP}-\inte C^{\PP}_n$. 
\end{lem}

\begin{proof} Assume that $T$ is in \mgp\ with respect to 
$\cup_{m\geq n}\bd C^{\PP}_m$. If the intersection is empty then $T$ lies in 
some $Y$ and hence bounds a compact \tm\ in $Y$. If the intersection is 
non-empty, then $T$ meets a single $\bd C^{\PP}_m$. Each annulus $A$ into 
which $T\cap\bd C^{\PP}_m$ splits $T$ must have $\bd A=\bd A\p$ for an 
annulus $A\p$ in $\bd C^{\PP}_m$. 

Consider an $A$ in $S=C^{\PP}_m-\inte C^{\PP}_n$. Let $T\p=A\cup A\p$. 
Then $T\p=\bd Q\p$ for a compact \tm\ $Q\p$ in $C^{\PP}_m$. We may assume 
that $Q\p\cap T=A$. Let $\widehat{S}$ and $\whC^{\PP}_m$ be obtained by 
adding a collar $C$ to these \tm s in $V^{\PP}-\inte C^{\PP}_m$. We may 
assme that $T$ meets $C$ in a product annulus. If $T\p$ is \inc\ in 
$\widehat{S}$, then $Q\p$ lies in $S$. If $T\p$ is compressible in 
$\widehat{S}$, then since $\widehat{S}$ is \irr\ $T\p$ bounds a solid torus 
or a classical knot space in $\widehat{S}$. This must be $Q\p$. So in either 
case $Q\p$ lies in $S$. Let $T\pp$ be the torus obtained from $T$ by 
replacing $A$ by $A\p$. Then $T\pp=\bd Q\pp$ for a compact \tm\ $Q\pp$ in 
$V^{\PP}$. If $T\pp$ is \inc\ in $V^{\PP}-\inte C^{\PP}_n$, 
then by induction $Q\pp$ lies in $V^{\PP}-\inte C^{\PP}_n$. If $T\pp$ is 
compressible in $V^{\PP}-\inte C^{\PP}_n$, then by irreducibility $T\pp$ 
bounds a solid torus or classical knot space in $V^{\PP}-\inte C^{\PP}_n$. 
This must be $Q\pp$. So in either case $Q\pp$ is in $V^{\PP}-\inte C^{\PP}_n$. 
If $Q\p\cap Q\pp=A\p$, then $T=\bd(Q\p\cup Q\pp)$. If $Q\p\cap Q\pp\neq A\p$, 
then $Q\p\sbs Q\pp$, and $T=\bd(Q\pp-\Inte Q\p)$. \end{proof}

%Proposition 7.8
\begin{prop} $V^{\PP}$ has finite genus. It has genus one if and only if 
\PP\ has exactly one element. \end{prop}

\begin{proof} $V^{\PP}$ has genus at most $m$. Since $V^p$ does not embed 
in \RRR\ the genus of $V^{\PP}$ must be at least one. So if $m=1$, then 
$V^p$ has genus one. Now suppose $m>1$. If $V^{\PP}$ has genus one, then it 
has a good exhaustion $\{K_n\}_{n\geq0}$ by solid tori. Choose $n$ and $k$ 
such that $K_0\sbs\inte C^{\PP}_n$ and $C^{\PP}_n\sbs\inte K_k$. Then since 
$\bd K_k$ is \inc\ in $V^{\PP}-\inte K_0$ it is \inc\ in the smaller space 
$V^{\PP}-\inte C^{\PP}_n$ and so bounds a compact \tm\ in this space, 
which is impossible. Thus $V^{\PP}$ has genus greater than one. \end{proof}

%Section 8
\section{The complex of \er s of $W$}

%Theorem 8.1
\begin{thm} Every $V^{\PP}$ is an \er\ of $W$ at each $C^{\PP}_n$. \end{thm}

\begin{proof} We know that $V^{\PP}$ is \ei\ rel $C^{\PP}_n$ in $W$. 
Clearly $W-V^{\PP}$ has no components with compact closure. Suppose $N$ is 
a regular \tm\ in $W$ such that $C^{\PP}_n\sbs\inte N$ and $\bd N$ is \inc\ in 
$W-C^{\PP}_n$. Then $N\sbs\inte C^{r,s}_m$ for some $r\leq p$, $s\geq q$, and 
$m>n$. We isotop $\bd N$ off a complete set of compressing disks for 
$\bd C^{r,s}_m$ in $C^{r,s}_m-\inte C^{\PP}_n$ so that $N$ lies in 
$C^{\PP}_m$. This can be done with compact support in $W-\inte C^{\PP}_n$. 
Running the isotopy backwards causes $V^{\PP}$ to engulf $N$. \end{proof}

%Theorem 8.2
\begin{thm} Let $V$ be an \er\ of $W$ at $J$, where $J\sbs\inte C^{\QQ}_n$. 
Then $V$ is isotopic to $V^{\PP}$ for some $\PP\sbs\QQ$ \end{thm}

\begin{proof} We may assume that $V$ is an \er\ of $W$ at a knot 
$\ka\sbs\inte J$. Let \PP\ be a minimal subset of \QQ\ such that, up to 
isotopy, $\ka\sbs\inte C^{\PP}_n$ for some $n$. Let \DD\ be the union of the 
set of co-cores of the 1-handles $H_{n,p}$ with $p\in\PP$. Then \ka\ is 
\DD-busting in $C^{\PP}_n$. 

If $m=1$, then clearly \ka\ is disk busting in $C^p_n$, so assume $m>1$. 

We let \whEE\ be the union of the attaching disks for the $H_{n+1,p}$ with 
$p\in\PP$ and the $D_{n+1,j}$ with $3p_1+1\leq j\leq 3p_m+1$. Let 
$\EE=\whEE\cap C^{\PP}_n$. We may assume that $\DD\sbs\EE$. The conditions 
of Proposition 5.3 are satisfied, so \ka\ is disk busting in $C^{\PP}_{n+1}$. 
It follows that $V$ is isotopic to $V^{\PP}$. \end{proof}

%Theorem 8.3
\begin{thm} $V^{\PP}$ and $V^{\QQ}$ are isotopic if and only if 
$\PP=\QQ$. \end{thm}

\begin{proof} We first consider the case $\PP=\{p\}$, $\QQ=\{q\}$, 
$p<q$. $V^p$ is an \er\ of $W$ at a knot \ka\ in $C^p_0$. Let \ta\ 
be the track of \ka\ under an isotopy taking $V^p$ to $V^q$ and 
\ka\ to $\ka\p$. Then $\ta\sbs\inte C^{r,s}_n$ for some $r\leq p$, 
$q\leq s$, and $n\geq0$. By the covering isotopy theorem 
\cite{Ch, EK} there is an ambient isotopy of \ka\ with 
track \ta\ which has compact support in $C^{r,s}_n$. 
Let $D$ be an attaching disk for $H_{n,p}$. Then \ka\ is $D$-busting 
in $C^{r,s}_n$, but $\ka\p$ is not. This is impossible since the isotopy 
is the identity on $\bd C^{r,s}_n$. 

Now consider the general case. Suppoes $p\in\PP$ and $p\notin\QQ$. Then 
$V^p$ is isotopic to $V^{\mathcal{R}}$ for some $\mathcal{R}\sbs\QQ$. 
Then we must have 
$\mathcal{R}=\{r\}$, where $r\neq p$, a contradiction. \end{proof}

%Theorem 8.4
\begin{thm} $V^{\PP}$ is minimal if and only if \PP\ has exactly one 
element. \end{thm}

\begin{proof} $V^p$ is clearly minimal. If $m>1$, then $V^{\PP}$ contains 
$V^p$ which is not \hm\ to $V^{\PP}$ since they have different genera. 
\end{proof}

%Theorem 8.5
\begin{thm} \SW\ is isomorphic to a triangulation of \R. \end{thm}

\begin{proof} The vertices of \SW\ are the $[V^p]$, $p\in\mathbf{Z}$. 
We have that $[V^p]$ and $[V^{p+1}]$ are joined by the edge $[V^{p,p+1}]$. 
Every \er\ of $W$ contained in $V^{p,p+1}$ is isotopic to $V^p$, $V^{p+1}$, 
or $V^{p,p+1}$. If $V$ is an \er\ of $W$ which contains representatives of 
$V^p$ and $V^q$, where $p<q$, then $V$ is isotopic to $V^{\PP}$, where 
$p,q\in\PP$. If $\PP\neq\{p,q\}$, then $V^{\PP}$ contains some $V^r$, 
$p\neq r\neq q$, so $[V]$ is not an edge joining $[V^p]$ and $[V^q]$. 
If $\PP=\{p,q\}$ and $q>p+1$, then \PP\ is bad, so $V$ is not \rirr, so 
again $[V]$ is not an edge. The result follows. \end{proof}

%Corollary 8.6
\begin{cor} If $W$ is a non-trivial covering space of a \tm\ $M$, then 
$\pi_1(M)\cong\mathbf{Z}$. \end{cor}

\begin{proof} This follows immediately from Theorem 8.5 and Corollary 3.4. 
\end{proof}

%Section 9
\section{Uncountably many $W$}

%Theorem 9.1
\begin{thm} There are uncountably many pairwise non-\hm\ $W$ each of which 
has all the properties of sections 7 and 8. \end{thm}

\begin{proof} Recall that all of the genus one \er s $V^p$ of a fixed 
$W$ resulting from our construction are \hm. We will modify our construction 
to obtain uncountably many $W$ such that different $W$ have non-\hm\ $V^p$. 

In our construction of $W\sh$ we used a copy of the same 2-tangle 
$\de^-\cup\de^+$ in $H$ for each 2-tangle $\de^-_n\cup\de^+_n$ in $H_{n+1}$. 
We will now change this so that the 2-tangle depends on $n$. 

We say that a \tm\ $Q$ is \textit{incompressibly embedded} in a \tm\ $X$ if 
$Q\sbs X$ and $\bd Q$ is \inc\ in $X$. 

%Lemma 9.2
\begin{lem} Given an excellent classical knot space $Q$, there is a 
\ps\ 2-tangle \ta\ in a 3-ball $B$ with exterior $X$ such that $Q$ is 
incompressibly embedded iin $X$ and every \inc\ torus in $X$ is isotopic 
to $\bd Q$. \end{lem}

\begin{proof} 
Let $B_0$ and $B_1$ be 3-balls. Let $D_i$ be a disk in $\bd B_i$. 
Let $A_i$ be an annulus in $\inte D_i$. Let $F_i$ be the annulus 
component of $D_i-\inte A_i$; let $E_i$ be the disk component. Let 
$\la^-_i\cup\la^+_i\cup\mu^-_i\cup\mu^+_i$ be a \pe\ 4-tangle in $B_i$. 
We require that $\la^{\pm}_0$ join $\bd B_0-\inte D_0$ to $\inte F_0$, 
$\mu^{\pm}_1$ join $\inte F_1$ to $\inte E_1$ , $\mu^{\pm}_0$ join $E_0$ to 
itself, and $\la^{\pm}_1$ join $\inte E_1$ to $\bd B_1-\inte D_1$. 
We now glue $B_0$ to $B_1$ by identifying $F_0\cup E_0$ with $F_1\cup E_1$ 
in such a way that $\la^{\pm}_0\cup\mu^{\pm}_1\cup\mu^{\pm}_0\cup\la^{\pm}_1$ 
is an arc $\de^{\pm}$. By Lemma 2.1 $\de^-\cup\de^+$ is a \pe\ system of two 
arcs in a 3-ball minus the interior of an unknotted solid torus with 
boundary $A_0\cup A_1$. We then glue $Q$ to this space by identifying $\bd Q$ 
with $A_0\cup A_1$ so that a meridian of $Q$ is glued to $\bd E_0=\bd E_1$. 
The result is a 3-ball $B$ cnotaining a 2-tangle $\ta=\de^-\cup\de^+$. 
Standard arguments then complete the proof. \end{proof}

Recall that $V^p$ is the monotone union of solid tori $C^p_n$, where 
$C^p_n=R_{n,p-1}\cup L_{n,p}\cup R_{n,p}\cup H_{n,p}$. Let 
$G_{n,p}=R_{n,p-1}\cup L_{n,p}\cup R_{n,p}$, 
$Y_{n+1,p}=C_{n+1,p}-\inte C_{n,p}$, $X_{n+1,p}=Y_{n+1,p}\cap H_{n+1,p}$, 
and $Z_{n+1,p}=Y_{n+1,p}\cap G_{n+1,p}$. 

Note that for all $n$ and $p$ the spaces $Z_{n+1,p}$ are \hm. 
It thus follows from the Jaco-Shalen-Johannson characteristic submanifold 
theory \cite{Ja, JS, Jo} that there are, up to 
homeomorphism, only finitely many excellent classical knot spaces which 
incompressibly embed in $Z_{n+1,p}$. Denote this set by \NN. 

Let \YY\ be the set of all homeomorphism types of excellent classical 
knot spaces which are not in \NN. For each infinite subset \SSS\ of \YY\ we 
construct a $W$ as follows. Choose a bijection of \SSS\ with the set 
$\mathbf{N}$ of natural numbers. For each $n\in\mathbf{N}$ use the 
corresponding knot space $Q_n$ in the construction of the 2-tangle $\ta_n$ 
in the previous lemma. Then use $\ta_n$ for $\de^-_{n-1}\cup\de^+_{n-1}$ 
in $H_n$. It follows that for each $n\geq m\geq0$ we have that $Q^p_{n+1}$ 
is incompressibly embedded in $V^p-\inte C^p_m$. 

%Lemma 9.3
\begin{lem} Suppose $Q\in\YY$ and $Q$ is incompressibly embedded in 
$V^p-\inte C^p_m$. Then $Q\in\SSS$. \end{lem}

\begin{proof} Since $Q$ is excellent it can be isotoped off 
$\cup_{n>m}\bd C^p_n$. It then lies in some $Y_{n+1,p}$. Since each 
$X_{n+1,p}$ is superb it can then be isotoped off $X_{n+1,p}\cap Z_{n+1,p}$. 
Since $Q\notin\NN$ it must lie in $X_{n+1,p}$ and thus be isotopic to 
$Q^p_{n+1}$. \end{proof}

Now suppose that $W\p$ is constructed using $\SSS\p$. Drop $p$ from the 
notation and denote the corresponding submanifolds of $W$ and $W\p$ by 
$V$ and $V\p$, $C_n$ and $C_n\p$, etc. 

%Lemma 9.4
\begin{lem} If $V$ and $V\p$ are \hm, then there are finite subsets 
$\SSS_0$ of \SSS\ and $\SSS_0\p$ of $\SSS\p$ such that 
$\SSS-\SSS_0=\SSS\p-\SSS_0\p$. \end{lem}

\begin{proof} Suppose $h:V\ra V\p$ is a homeomorphism. Choose $m$ and $k$ 
such that $h(C_0)\sbs\inte C\p_m$ and $C\p_m\sbs\inte h(C_k)$. Then for all 
$n\geq k$ we have that $h(\bd C_n)$ is \inc\ in $V\p-\inte h(C_0)$ and 
hence is \inc\ in the smaller space $V\p-\inte C\p_m$. It follows that 
$h(Q_{n+1})$ is isotopic in this space to some $Q\p_{j+1}$ with $j\geq m$. 
Let $\AAA=\{Q_1, \ldots, Q_k\}$. Then $\SSS-\AAA\sbs\SSS\p$. A similar 
argument using $h\n$ yields a finite set $\AAA\p\sbs\SSS\p$ such that 
$\SSS\p-\AAA\p\sbs\SSS$. We then let $\SSS_0=\AAA\cup(\SSS\cap\AAA\p)$ 
and $\SSS\p_0=\AAA\p\cup(\SSS\p\cap\AAA)$. \end{proof}

Define an equivalence relation on the 
set of infinite subsets of \YY\ by setting $\SSS\sim\SSS\p$ if 
$\SSS-\SSS_0=\SSS\p-\SSS\p_0$ as in the lemma. Each equivalence class has only 
countably many elements, and so there are uncountably many equivalence 
classes. It follows that there are uncountably many non-\hm\ $V$ 
and hence uncountably many non-\hm\ $W$. \end{proof}

\end{document}